\documentclass[fleqn]{article}
\usepackage{amsfonts}
\usepackage{amsmath}

\newtheorem{theorem}{Theorem}

\newtheorem{corollary}[theorem]{Corollary}

\newtheorem{definition}[theorem]{Definition}
\newtheorem{example}[theorem]{Example}

\newenvironment{proof}[1][Proof]{\noindent\textbf{#1.} }{\ \rule{0.5em}{0.5em}}
\input{tcilatex}
\begin{document}

\title{Osculating mate of a Frenet curve in the Euclidean 3-space}
\author{Ak{\i}n Alkan$^1$, Mehmet \"{O}nder$^2$ \\
%EndAName
$^1$Manisa Celal Bayar University, \\
G\"{o}rdes Vocational School, 45750, G\"{o}rdes, Manisa, Turkey.\\
$^2$Delibekirli Village, Tepe Street, 31440, K{\i}r{\i}khan, Hatay, Turkey.\\
E-mails: $^1$akin.alkan@cbu.edu.tr, $^2$mehmetonder197999@gmail.com \\
Orcid Ids: $^1$https://orcid.org/0000-0002-8179-9525, \\
$^2$https://orcid.org/0000-0002-9354-5530}
\maketitle

\begin{abstract}
A new kind of partner curve called osculating mate of a Frenet curve is
introduced. Some characterizations for osculating mate are obtained and
using the obtained results some special curves such as slant helix,
spherical helix, $C$-slant helix and rectifying curve are constructed.
\end{abstract}

\textbf{AMS Classsification:} 53A04, 53C40.


\bigskip

\section{Introduction}

The most fascinating and important subject of curve theory is to obtain the
characterizations for a curve or a curve pair which are known as special
curves or partner curves. Helices, slant helices, rectifying curves,
spherical curves, etc. are common examples of such curves. Especially, the
helices are seen in many areas such as nature, design of mechanic tools and
highways, simulation of kinematic motion or architect, nucleic acids and
molecular model of DNA \cite{Peyrard,Rapaport,Walsby1,Walsby2,Watson}.
Helices are also important in physics since they are used in helical gears,
shapes of springs and elastic rods \cite{Healey,Keil}. A helix $\alpha $ is
defined by that the tangent of $\alpha $ always makes a constant angle with
a fixed direction and necessary and sufficient condition for a curve $\alpha 
$ to be a helix is that $\frac{\tau }{\kappa }(s)$ is constant, where $%
\kappa $ is the first curvature (or curvature) and $\tau $ is the second
curvature (or torsion) of $\alpha $ \cite{Barros,Struik}. Another kind of
special curves is slant helix defined by that there exists always a constant
angle between the principal normal line of curve and a fixed direction. This
special curve was first defined by Izumiya and Takeuchi \cite{Izumiya1}.
Later, Z\i plar, \c{S}enol and Yayl\i\ have introduced a new special curve
called Darboux helix and they have obtained that a curve is a Darboux helix
iff the curve is a slant helix \cite{Ziplar}.

Furthermore, a special curve can be defined by considering its position
vector. A curve $\alpha $ in the Euclidean 3-space $E^{3}$ for which the
position vector of the curve is always contained in its rectifying plane
(respectively, osculating plane or normal plane) is named as rectifying
curve or briefly rectifying (respectively, osculating curve or normal curve) 
\cite{Chen1}. Rectifying curves, normal curves and osculating curves satisfy
Cesaro's fixed point condition \cite{Otsuki}. Namely, rectifying, normal and
osculating planes of such curves always contain a particular point.
Moreover, Darboux vectors (centrodes) and rectifying curves are related and
used in different areas of sciences such as kinematics, mechanics and
differential geometry of curves of constant precession \cite{Chen2}.

K\i z\i ltu\u{g}, \"{O}nder and Yayl\i\ have defined a new kind of special
curves called normal direction curves\cite{Kiziltug}. Later, \c{C}akmak has
studied the same subject and two similar ones in 3-dimensional compact Lie
group \cite{AliCakmak}.

Recently, Deshmukh, Chen and Alghanemi have studied natural mate and
conjugate mate of a curve \cite{Deschmukh}. They have given some new
characterizations for spherical curve, helix, rectifying curve and slant
helix. Alghanemi and Khan have given the position vectors of natural mate
and conjugate mate \cite{Algha}. Mak has studied these mates in
three-dimensional Lie groups \cite{Mak}. Later, Camc\i\ et all have studied
sequential natural mates of Frenet curves in $E^{3}$\cite{Camci}.

In the present paper, we define osculating mate of a Frenet curve $\alpha $
in $E^{3}$. We give some relations between a Frenet curve and its osculating
mate and introduce some applications of osculating mates to slant helix,
spherical helix, rectifying curve and $C$-slant helix in $E^{3}$.

\bigskip

\section{Preliminaries}

Let $\alpha :I\rightarrow E^{3}~$be a unit speed curve with arclength
parameter $s$. The vector $T(s)=\alpha ^{\prime }(s)$ is called unit tangent
vector of $\alpha ~$and the function $\kappa (s)=\left\Vert \alpha ^{\prime
\prime }(s)\right\Vert $ is called the curvature of $\alpha .$ The unit
principal normal vector $N(s)~$of the curve $\alpha ~$is defined by $\alpha
^{\prime \prime }(s)=\kappa (s)N(s)$. The unit binormal vector of $\alpha $
is $B(s)=T(s)\times N(s).$Then, the Frenet frame $\left\{ T,N,B\right\} ~$%
has the following formulas 
\begin{equation*}
\left( 
\begin{array}{c}
T^{\prime } \\ 
N^{\prime } \\ 
B^{\prime }%
\end{array}%
\right) =\left( 
\begin{array}{ccc}
0 & \kappa & 0 \\ 
-\kappa & 0 & \tau \\ 
0 & -\tau & 0%
\end{array}%
\right) \left( 
\begin{array}{c}
T \\ 
N \\ 
B%
\end{array}%
\right)
\end{equation*}%
where $\tau =\tau (s)~$is the torsion of the curve $\alpha $\ and defined by 
$\tau =-\left\langle B^{\prime },N\right\rangle $\cite{Struik}. If $\kappa
(s)\neq 0$, the curve $\alpha ~$is named as Frenet curve. The curve $\alpha
~ $is a general helix iff $\frac{\tau }{\kappa }(s)$ is constant. Similarly,
the characterization of a slant helix is given by the necessary and
sufficient condition that 
\begin{equation}
\sigma (s)=\frac{\kappa ^{2}}{(\kappa ^{2}+\tau ^{2})^{3/2}}\left( \frac{%
\tau }{\kappa }\right) ^{\prime }=const.  \label{eqn1}
\end{equation}%
(\cite{Izumiya1}).

A Frenet curve $\alpha ~$is named as a Salkowski(respectively,
anti-Salkowski) curve if its curvature $\kappa $ is constant but torsion $%
\tau $ is non-constant(respectively, torsion $\tau $ is constant but
curvature $\kappa $ is non-constant) \cite{Monterde2}.

A Frenet curve $\alpha ~$is named as a spherical curve if all points of $%
\alpha $ lie on the same sphere and such a curve is characterized as
follows: A Frenet curve $\alpha $ is a spheciral curve iff $(p^{\prime
}q)^{\prime }+\frac{p}{q}=0$ holds, where $p=1/\kappa ,~q=1/\tau $.
Moreover, another characterization for a spherical curve is that a Frenet
curve $\alpha $ is a spheciral curve iff $p^{2}+(p^{\prime }q)^{2}=a^{2}$
holds, where $a>0$ is the radius of the spehere on which $\alpha $ lies \cite%
{Millman}.

A Frenet curve $\alpha ~$is named as rectifying curve if the position vector
of $\alpha $ always lies on the rectifying plane of the curve \cite{Chen1}.
A rectifying curve is characterized by the necessary and sufficiant
condition that $\frac{\tau }{\kappa }(s)=\frac{1}{c}(s+b)$ holds, where $%
c\neq 0,~b$ are real constants and such a curve has the parametrization $%
\alpha (s)=(s+b)T(s)+cB(s)$ \cite{Chen1}.

The vector $W$ defined by $W=\frac{\tau T+\kappa B}{\sqrt{\kappa ^{2}+\tau
^{2}}}$ is called unit Darboux vector of $\alpha $. Then, the frame $\left\{
N,C=W\times N,W\right\} $ is called the alternative frame of $\alpha .$ A
curve $\alpha $ is called Darboux helix if the unit Darboux vector $W$ makes
a constant angle with a fixed direction and the curve $\alpha $ is a Darboux
helix iff $\alpha $ is a slant helix \cite{Ziplar}. A curve $\alpha $ is
named as $C$-slant helix if the unit vector $C$ always makes a constant
angle with a fixed direction. Necessary and sufficient condition for a curve 
$\alpha $ to be a $C$-slant helix is that the function 
\begin{equation}
\mu (s)=\frac{(f^{2}+g^{2})^{3/2}}{f^{2}(\frac{g}{f})^{\prime }},
\label{eqn1-1}
\end{equation}%
is constant \cite{Uzunoglu}.

\qquad

\section{Osculating mates of a Frenet curve in $E^{3}$}

In this section, we define osculating mate of a Frenet curve in $E^{3}$ and
give some characterizations for this curve.

\begin{definition}
Let $\alpha :I\subset 
%TCIMACRO{\U{211d} }%
%BeginExpansion
\mathbb{R}
%EndExpansion
\rightarrow E^{3}~$be a unit speed Frenet curve. The curve $\beta ~$defined
by%
\begin{equation}
\beta (s)=\int \left( x_{1}(s)T(s)+x_{2}(s)N(s)\right) ,  \label{eqn2}
\end{equation}%
and satisfying the conditions $x_{1}^{2}(s)+x_{2}^{2}(s)=1~$and $\beta
^{\prime \prime }\perp sp\left\{ T,N\right\} ~$is named as the osculating
mate of the curve $\alpha $.
\end{definition}

\bigskip Unless otherwise stated, hereafter when we talk about the concept
of curves we will mean Frenet curves.

\begin{theorem}
The Frenet apparatus of osculating mate $\beta $ are computed as follows%
\begin{equation}
\left\{ 
\begin{array}{c}
\overset{\_}{T}=\sin \left( \int \kappa (s)ds\right) T+\cos \left( \int
\kappa (s)ds\right) N,~\overset{\_}{N}=B,~ \\ 
\overset{\_}{B}=\cos \left( \int \kappa (s)ds\right) T-\sin \left( \int
\kappa (s)ds\right) N,%
\end{array}%
\right.  \label{eqn3}
\end{equation}%
\begin{equation}
\overset{\_}{\kappa }=\varepsilon _{1}\tau \cos \left( \int \kappa
(s)ds\right) ,~\overset{\_}{\tau }=\tau \sin \left( \int \kappa (s)ds\right)
,  \label{eqn4}
\end{equation}%
where $\varepsilon _{1}=\pm 1$ is chosen such as $\overset{\_}{\kappa }>0.$
\end{theorem}

\begin{proof}
Let the Frenet apparatus of osculating mate $\beta $ be given by $\left\{ 
\overset{\_}{T},\overset{\_}{N},\overset{\_}{B};\overset{\_}{\kappa },%
\overset{\_}{\tau }\right\} $. From Definition 1, it follows $\beta ^{\prime
}=\overset{\_}{T}=x_{1}T+x_{2}N$ . Differentiating last equality we have,%
\begin{equation}
\overset{\_}{T}^{\prime }=\left( x_{1}^{\prime }-x_{2}\kappa \right)
T+\left( x_{2}^{\prime }+x_{1}\kappa \right) N+x_{2}\tau B,  \label{eqn5}
\end{equation}%
which gives the system%
\begin{equation}
x_{1}^{\prime }-x_{2}\kappa =0,~x_{2}^{\prime }+x_{1}\kappa =0,~x_{2}\tau
\neq 0.  \label{eqn6}
\end{equation}%
The solution of the system (\ref{eqn6}) is%
\begin{equation}
x_{1}(s)=\sin \left( \int \kappa (s)ds\right) ,~x_{2}(s)=\cos \left( \int
\kappa (s)ds\right) .  \label{eqn7}
\end{equation}%
Then, it follows $\overset{\_}{T}=\sin \left( \int \kappa ds\right) T+\cos
\left( \int \kappa ds\right) N~$and from (\ref{eqn5}), we have $\overset{\_}{%
\kappa }\overset{\_}{N}=\tau \cos \left( \int \kappa ds\right) B$. Hence, we
obtain%
\begin{equation}
\overset{\_}{\kappa }=\varepsilon _{1}\tau \cos \left( \int \kappa ds\right)
,~\overset{\_}{N}=B,  \label{eqn8}
\end{equation}%
where $\varepsilon _{1}=\pm 1$ is chosen such as $\overset{\_}{\kappa }>0.~$%
Furthermore,%
\begin{equation}
\overset{\_}{B}=\overset{\_}{T}\times \overset{\_}{N}=\cos \left( \int
\kappa ds\right) T-\sin \left( \int \kappa ds\right) N.  \label{eqn9}
\end{equation}%
Differentiating (\ref{eqn9}) and using the equality $\overset{\_}{\tau }%
^{\prime }=-\left\langle \overset{\_}{B}^{\prime },\overset{\_}{N}%
\right\rangle $, we have $\overset{\_}{\tau }=\tau \sin \left( \int \kappa
ds\right) .$
\end{proof}

\bigskip

\begin{theorem}
The curvatures $\kappa $ and $\tau $ of $\alpha $ are computed as%
\begin{equation}
\kappa =\frac{\varepsilon _{1}\overset{\_}{\kappa }^{2}}{\overset{\_}{\kappa 
}^{2}+\overset{\_}{\tau }^{2}}\left( \frac{\overset{\_}{\tau }}{\overset{\_}{%
\kappa }}\right) ^{\prime },~\tau =\pm \sqrt{\overset{\_}{\kappa }^{2}+%
\overset{\_}{\tau }^{2}},  \label{eqn10}
\end{equation}%
respectively.
\end{theorem}

\begin{proof}
From (\ref{eqn4}), we easily get%
\begin{equation}
\tau =\pm \sqrt{\overset{\_}{\kappa }^{2}+\overset{\_}{\tau }^{2}}.
\label{eqn11}
\end{equation}%
Writing (\ref{eqn11}) into (\ref{eqn4}), it follows%
\begin{equation}
\cos \left( \int \kappa ds\right) =\frac{\pm \varepsilon _{1}\overset{\_}{%
\kappa }}{\sqrt{\overset{\_}{\kappa }^{2}+\overset{\_}{\tau }^{2}}},~\sin
\left( \dint \kappa ds\right) =\frac{\pm \overset{\_}{\tau }}{\sqrt{\overset{%
\_}{\kappa }^{2}+\overset{\_}{\tau }^{2}}},  \label{eqn12}
\end{equation}%
respectively. By taking the derivative of the second equality in (\ref{eqn12}%
), we get%
\begin{equation}
\kappa \cos \left( \int \kappa ds\right) =\pm \frac{\overset{\_}{\kappa }(%
\overset{\_}{\kappa }\overset{\_}{\tau }^{\prime }-\overset{\_}{\kappa }%
^{\prime }\overset{\_}{\tau })}{(\overset{\_}{\kappa }^{2}+\overset{\_}{\tau 
}^{2})^{3/2}}.  \label{eqn13}
\end{equation}%
Writing first equality in (\ref{eqn12}) into (\ref{eqn13}) gives $\kappa =%
\frac{\varepsilon _{1}\overset{\_}{\kappa }^{2}}{\overset{\_}{\kappa }^{2}+%
\overset{\_}{\tau }^{2}}\left( \frac{\overset{\_}{\tau }}{\overset{\_}{%
\kappa }}\right) ^{\prime },~$which completes the proof.
\end{proof}

\bigskip

From Theorem 3, Theorem 4 and equation (\ref{eqn1}), we have%
\begin{equation}
\overset{\_}{\tau }=\tau \sin \left( \int \kappa ds\right) ,~\frac{\kappa }{%
\tau }=\varepsilon _{1}\overset{\_}{\sigma },  \label{eqn14}
\end{equation}%
which gives the following corollary.

\begin{corollary}
i) $\alpha $\ is plane curve iff the osculating mate $\beta $ is plane curve.

ii) $\alpha $\ is helix iff the osculating mate $\beta $ is slant helix.
\end{corollary}

\begin{theorem}
The osculating mate $\beta $ is spherical curve iff the curvatures $\kappa
,~\tau $ of $\alpha $ satisfy the following equality%
\begin{equation}
\left( \tau \cos x\right) ^{\prime }=\pm \tau ^{2}\sin x\cos x\sqrt{%
a^{2}\tau ^{2}\cos ^{2}x-1},  \label{eqn15}
\end{equation}%
where $a>0$ is the radius of the sphere and $x(s)=\dint \kappa (s)ds$.
\end{theorem}

\begin{proof}
First assume that $\beta $ lies on a sphere with radius $a>0$. Hence, $%
\overset{\_}{p}^{2}+(\overset{\_}{p}^{\prime }\overset{\_}{q})^{2}=a^{2},~$%
where $\overset{\_}{p}=1/\overset{\_}{\kappa },~\overset{\_}{q}=1/\overset{\_%
}{\tau }$. From (\ref{eqn4}), it follows $\overset{\_}{p}^{\prime }=\frac{%
-\varepsilon _{1}(\tau \cos x)^{\prime }}{\tau ^{2}\cos ^{2}x}$. Hence, we
have 
\begin{equation}
\frac{1}{\tau ^{2}\cos ^{2}x}\left[ 1+\frac{((\tau \cos x)^{\prime })^{2}}{%
\tau ^{4}\sin ^{2}x\cos ^{2}x}\right] =a^{2},  \label{eqn16}
\end{equation}%
which gives (\ref{eqn15}).

Conversely, assume that (\ref{eqn15}) holds. By differentiating the first
equality in (\ref{eqn4}), we have 
\begin{equation}
-\overset{\_}{p}^{\prime }=\frac{\overset{\_}{\kappa }^{\prime }}{\overset{\_%
}{\kappa }^{2}}=\frac{-\varepsilon _{1}(\tau \cos x)^{\prime }}{\tau
^{2}\cos ^{2}x}.  \label{eqn17}
\end{equation}%
Writing (\ref{eqn15}) in (\ref{eqn17}) gives%
\begin{equation}
\overset{\_}{p}^{\prime }=\frac{\mp \varepsilon _{1}\sin x\sqrt{a^{2}\tau
^{2}\cos ^{2}x-1}}{\cos x}.  \label{eqn18}
\end{equation}%
Then, by taking into account the second equality in (\ref{eqn4}), we obtain $%
\overset{\_}{p}^{\prime }\overset{\_}{q}=\frac{\mp \varepsilon _{1}\sqrt{%
a^{2}\tau ^{2}\cos ^{2}x-1}}{\tau \cos x}$ and so, $\overset{\_}{p}^{2}+(%
\overset{\_}{p}^{\prime }\overset{\_}{q})^{2}=a^{2}$, i.e., $\beta $ lies on
a sphere with radius $a>0$.
\end{proof}

\begin{theorem}
The osculating mate $\beta $ is rectifying iff the function $\tan \int
\kappa ds$ is a linear function of $s.$
\end{theorem}

\begin{proof}
Suppose that $\beta $ is rectifying. So, we have $\frac{\overset{\_}{\tau }}{%
\overset{\_}{\kappa }}=\frac{1}{c}\left( s+b\right) ,$ where $c\neq 0,~b$
are real constants. Considering (\ref{eqn4}), it follows $\tan \int \kappa
ds=\frac{\varepsilon _{1}}{c}\left( s+b\right) $.

Conversely, let we write $\varepsilon _{1}\tan \int \kappa ds=\left(
a_{1}s+a_{2}\right) ,$ where $a_{1}\neq 0,~a_{2\text{ }}$are real constants.
Let define $a_{1}=\frac{1}{c}$ and $a_{2}=\frac{b}{c},$ where $c\neq 0$ is a
real contant. Then, we get $\varepsilon _{1}\tan \int \kappa ds=\frac{1}{c}%
\left( s+b\right) $ and it follows $c\tau \sin \left( \int \kappa ds\right)
=\varepsilon _{1}\left( s+b\right) \tau \cos \left( \int kds\right) $. By
taking into account (\ref{eqn4})$,$ we obtain $\left( s+b\right) \overset{\_}%
{\kappa }-c\overset{\_}{\tau }=0,$ which gives that $\beta $ is rectifying.
\end{proof}

\begin{theorem}
The position vector of osculating mate $\beta $ is given by 
\begin{equation}
\beta =\left[ \int \left( -\frac{\kappa }{\tau }h^{\prime }+\sin \left( \int
\kappa ds\right) \right) ds\right] T-\frac{h^{\prime }}{\tau }N+hB,
\label{eqn19}
\end{equation}%
where $h(s)=\frac{(dd^{\prime })^{\prime }-1}{\tau \cos \int \kappa ds}$ and 
$d=d(s)=\left\Vert \beta (s)\right\Vert $ is the distance function of $\beta
.$
\end{theorem}

\begin{proof}
For the position vector $\beta ,$ we can write 
\begin{equation}
\beta =a_{1}T+a_{2}N+a_{3}B,  \label{eqn20}
\end{equation}%
where $a_{i}=a_{i}(s),~\left( i=1,2,3\right) ~$are smooth functions of $s$.
Differentiating (\ref{eqn20}) and using (\ref{eqn3}), we have 
\begin{equation}
\left\{ 
\begin{array}{c}
\sin \left( \int \kappa ds\right) T+\cos \left( \int \kappa ds\right)
N=\left( a_{1}^{\prime }-a_{2}\kappa \right) T+\left( a_{1}\kappa
+a_{2}^{\prime }-a_{3}\tau \right) N \\ 
+\left( a_{2}\tau +a_{3}^{\prime }\right) B.%
\end{array}%
\right.  \label{eqn21}
\end{equation}%
From (\ref{eqn21}), we have the following system 
\begin{equation}
\left\{ 
\begin{array}{c}
a_{1}^{\prime }-a_{2}\kappa =\sin \left( \int \kappa ds\right) , \\ 
a_{1}\kappa +a_{2}^{\prime }-a_{3}\tau =\cos \left( \int \kappa ds\right) ,
\\ 
a_{2}\tau +a_{3}^{\prime }=0.%
\end{array}%
\right.  \label{eqn22}
\end{equation}%
From (\ref{eqn20}), it follows $d^{2}=a_{1}^{2}+a_{2}^{2}+a_{3}^{2}.$
Differentiating last equality gives $dd^{\prime }=$ $a_{1}a_{1}^{\prime
}+a_{2}a_{2}^{\prime }+a_{3}a_{3}^{\prime }$. Then, from system (\ref{eqn22}%
), we get 
\begin{equation}
dd^{\prime }=a_{1}\sin \left( \int \kappa ds\right) +a_{2}\cos \left( \int
\kappa ds\right) .  \label{eqn23}
\end{equation}%
Differentiating (\ref{eqn23}) and taking into account system (\ref{eqn22}),
we obtain $a_{3}=\frac{(dd^{\prime })^{\prime }-1}{\tau \cos \left( \int
\kappa ds\right) }.$ By writing $h(s)=a_{3}(s),$ from system (\ref{eqn22}),
we get%
\begin{equation}
a_{2}=-\frac{h^{\prime }}{\tau }(dd^{\prime })^{\prime },a_{1}=\int \left[ -%
\frac{\kappa }{\tau }\left( h\right) ^{\prime }+\sin \left( \int \kappa
ds\right) \right] ds.  \label{eqn24}
\end{equation}%
Considering (\ref{eqn20}), we have (\ref{eqn19}).
\end{proof}

\begin{corollary}
Let $\beta $ be an osculating mate of $\alpha .$

i) $\beta $ is spherical curve iff $h(s)=\frac{-1}{\tau \cos \int \kappa ds}%
. $

ii) If $\beta $ is rectifying curve, then $h=0.$
\end{corollary}

\begin{proof}
i) $\beta $ is spherical curve iff $d$ is a non-zero constant iff $h(s)=%
\frac{-1}{\tau \cos \int \kappa ds}.$

ii) Since $\beta $ is a rectifying curve, its distance function $d$
satisfies $d^{2}(s)=s^{2}+c_{1}s+c_{2},$ where $c_{i};(i=1,2)$ are constants 
\cite{Chen2}. Then, we have $h=0$.
\end{proof}

\begin{theorem}
\bigskip Let $\beta $ be an osculating mate of $\alpha .$\ 

i) $\beta $ is Bertrand curve iff the function $(pq^{\prime })^{2}+q^{2}$ is
a non-zero constant.

ii) $\alpha $ is Bertrand curve iff $\varepsilon _{1}\varsigma _{1}\overset{%
\_}{\sigma }\mp \varsigma _{2}=\frac{1}{\sqrt{\overset{\_}{\kappa }^{2}+%
\overset{\_}{\tau }^{2}}}$, where $\varsigma _{1}\neq 0,$ $\varsigma _{2}$
are constants.
\end{theorem}

\begin{proof}
i) Since $\beta $ is a Bertrand curve, we can write $a\overset{\_}{\kappa }+b%
\overset{\_}{\tau }=1,$ where $a\neq 0$ and $b$ are constants \cite{Bertrand}%
. Writing (\ref{eqn4}) in the last equality gives 
\begin{equation}
a\varepsilon _{1}\cos \left( \int \kappa ds\right) +b\sin \left( \dint
\kappa ds\right) =\frac{1}{\tau }=q.  \label{eqn24-1}
\end{equation}%
By differentiating (\ref{eqn24-1}), we have 
\begin{equation}
-a\varepsilon _{1}\sin \left( \int \kappa ds\right) +b\cos \left( \dint
\kappa ds\right) =\left( \frac{1}{\tau }\right) ^{\prime }\frac{1}{\kappa }%
=q^{\prime }p.  \label{eqn24-2}
\end{equation}%
From (\ref{eqn24-1}) and (\ref{eqn24-2}), it follows $(pq^{\prime
})^{2}+q^{2}=a^{2}+b^{2}.$

Conversely, let $(pq^{\prime })^{2}+q^{2}$ be a non-zero constant. Define $%
q= $ $(a^{2}+b^{2})\cos \theta $ and $pq^{\prime }=(a^{2}+b^{2})\sin \theta
, $ where $a\neq 0,~b$ are real constants. Differentiating first equality
and writing the result in the second one gives $\theta ^{\prime }=-\kappa $.
Then, the equality $q=$ $(a^{2}+b^{2})\cos \theta ~$becomes $\tau =\frac{1}{%
(a^{2}+b^{2})\cos \left( \dint \kappa ds+m\right) }$, where $m$ is
integration constant. By taking into account (\ref{eqn4}), we have 
\begin{equation}
\overset{\_}{\kappa }=\frac{\varepsilon _{1}\cos \left( \dint \kappa
ds\right) }{\left( a^{2}+b^{2}\right) \cos \left( \dint \kappa ds+m\right) }%
,~\overset{\_}{\tau }=-\frac{\sin \left( \dint \kappa ds\right) }{\left(
a^{2}+b^{2}\right) \cos \left( \dint \kappa ds+m\right) }.  \label{eqn24-3}
\end{equation}%
By writing $A=\left( a^{2}+b^{2}\right) \cos (m)$, $B=\left(
a^{2}+b^{2}\right) \sin (m)$ and taking into account (\ref{eqn24-3}) it
follows $A\overset{\_}{\kappa }+B\overset{\_}{\tau }=1$, i.e., $\beta $ is
Bertrand curve.

ii) If $\alpha $ is a Bertrand curve, then $\varsigma _{1}\kappa +\varsigma
_{2}\tau =1,$where $\varsigma _{1}\neq 0,$ $\varsigma _{2}$ are constants.
Writing (\ref{eqn10}) in the last equality, it follows $\frac{\varepsilon
_{1}\varsigma _{1}\overset{\_}{\kappa }^{2}}{\overset{\_}{\kappa }^{2}+%
\overset{\_}{\tau }^{2}}\left( \frac{\overset{\_}{\tau }}{\overset{\_}{%
\kappa }}\right) ^{\prime }\pm \varsigma _{2}\sqrt{\overset{\_}{\kappa }^{2}+%
\overset{\_}{\tau }^{2}}=1$ or equivalently, $\varepsilon _{1}\varsigma _{1}%
\overset{\_}{\sigma }\pm \varsigma _{2}=\frac{1}{\sqrt{\overset{\_}{\kappa }%
^{2}+\overset{\_}{\tau }^{2}}}.$

Conversely, if $\varepsilon _{1}\varsigma _{1}\overset{\_}{\sigma }\pm
\varsigma _{2}=\frac{1}{\sqrt{\overset{\_}{\kappa }^{2}+\overset{\_}{\tau }%
^{2}}}$ holds, by taking into account (\ref{eqn10}), we have $\varsigma
_{1}\kappa +\varsigma _{2}\tau =1$, i.e., $\alpha $ is a Bertrand curve.
\end{proof}

\begin{theorem}
Let $\beta $ be an osculating mate of $\alpha .$

i) $\beta $ is Mannheim curve iff $\frac{1}{\tau }\cos \int \kappa
ds=\varepsilon _{1}\lambda _{1},$ where $\lambda _{1}$ is non-zero constant.

ii) $\alpha $ is Mannheim curve iff $\varepsilon _{1}\sqrt{\overset{\_}{%
\kappa }^{2}+\overset{\_}{\tau }^{2}}\overset{\_}{\sigma }^{3}=\lambda
_{2}\left( 1+\overset{\_}{\sigma }^{2}\right) ,$ where $\lambda _{2}$ is
non-zero constant.
\end{theorem}

\begin{proof}
i) If $\beta $ is Mannheim curve, there is a non-zero constant $\lambda _{1}$
such that $\overset{\_}{\kappa }=\lambda _{1}(\overset{\_}{\kappa }^{2}+%
\overset{\_}{\tau }^{2})$ holds \cite{Liu,Honda}. Writing (\ref{eqn4}) in
the last equality gives $\frac{1}{\tau }\cos \int \kappa ds=\varepsilon
_{1}\lambda _{1}.$

Conversely, if $\frac{1}{\tau }\cos \int \kappa ds=\varepsilon _{1}\lambda
_{1}$ holds for a non-zero constant $\lambda _{1},$ from (\ref{eqn3}) and (%
\ref{eqn4}), we have that $\overset{\_}{\kappa }=\lambda _{1}(\overset{\_}{%
\kappa }^{2}+\overset{\_}{\tau }^{2})$ holds, i.e., $\beta $ is Mannheim
curve.

ii) If $\alpha $ is a Mannheim curve, the curvatures of $\alpha $ satisfy $%
\kappa =\lambda _{2}(\kappa ^{2}+\tau ^{2}),$ where $\lambda _{2}$ is
non-zero constant. Hence, we get $\frac{1}{\kappa }=\lambda _{2}(1+\frac{%
\tau ^{2}}{\kappa ^{2}}).$ Writing (\ref{eqn10}) in the last equality gives
and considering (\ref{eqn1}), we obtain $\varepsilon _{1}\sqrt{\overset{\_}{%
\kappa }^{2}+\overset{\_}{\tau }^{2}}\overset{\_}{\sigma }^{3}=\lambda
_{2}\left( 1+\overset{\_}{\sigma }^{2}\right) .$

The converse is clear. \bigskip
\end{proof}

\begin{corollary}
Let $\beta $ be an osculating mate of $\alpha .$Then, $\alpha $ is Mannheim
curve iff he cuvatures of $\alpha $ and $\beta $ satisfy $\overset{\_}{%
\kappa }=\pm \lambda \tau ,$ where $\lambda $ is non-zero constant.
\end{corollary}

\begin{theorem}
\bigskip The curve $\beta $ be an osculating mate of $\alpha .$

i) Let $\alpha $ be Salkowski curve. Then $\beta $ is Salkowski curve iff $%
\tau =\varepsilon _{1}e_{3}\sec (e_{1}s+e_{2})$ , where $e_{i};~(i=1,2,3)$
are real constants.

ii) Let $\beta $ be Salkowski curve with constant curvature $\overset{\_}{%
\kappa }=e_{4}$.Then $\alpha $ is Salkowski curve with $\kappa =c>0$ iff $%
\varepsilon _{1}e_{4}\overset{\_}{\tau }^{\prime \prime }-2c\overset{\_}{%
\tau }\overset{\_}{\tau }^{\prime }=0$ holds.

iii) Let $\beta $ be anti-Salkowski curve with constant torsion $\overset{\_}%
{\tau }=e_{5}$.Then $\alpha $ is Salkowski curve with $\kappa =c>0$ iff $%
\varepsilon _{1}e_{5}\overset{\_}{\kappa }^{\prime \prime }+2c\overset{\_}{%
\kappa }\overset{\_}{\kappa }^{\prime }=0$ holds.
\end{theorem}

\begin{proof}
i) Since $\alpha $ is Salkowski curve, we have $\kappa =e_{1}>0$ is constant
but $\tau $ is non-constant. Then, from (\ref{eqn3}) it follows $\overset{\_}%
{\kappa }=\varepsilon _{1}\tau \cos \left( e_{1}s+e_{2}\right) ,~\overset{\_}%
{\tau }=\tau \sin \left( e_{1}s+e_{2}\right) ,~$where $e_{2}$ is integration
contant. So, we get $\frac{\overset{\_}{\tau }}{\overset{\_}{\kappa }}%
=\varepsilon _{1}\tan \left( e_{1}s+e_{2}\right) $. Hence, $\beta $ is
Salkowski curve with constant curvature $\overset{\_}{\kappa }=e_{3}>0$ iff $%
\tau =\varepsilon _{1}e_{3}\sec (e_{1}s+e_{2}).$

The proofs of (ii) and (iii) are similar to proof of (i).
\end{proof}

\bigskip

Let now $(\overset{\_}{T}),~(\overset{\_}{N}),~(\overset{\_}{B})$ denote the
tangent indicatrix, the principal normal indicatrix and the binormal
indicatrix of osculating mate $\beta ,$ respectively. Then, the curvatures
of these spherical curves are computed as 
\begin{equation}
\kappa _{\overset{\_}{T}}=\frac{\sqrt{\overset{\_}{\kappa }^{2}+\overset{\_}{%
\tau }^{2}}}{\overset{\_}{\kappa }},~\tau _{\overset{\_}{T}}=\frac{\overset{%
\_}{\kappa }}{\overset{\_}{\kappa }^{2}+\overset{\_}{\tau }^{2}}\left( \frac{%
\overset{\_}{\tau }}{\overset{\_}{\kappa }}\right) ^{\prime },  \label{eqn25}
\end{equation}

\begin{equation}
\kappa _{\overset{\_}{N}}=\frac{\sqrt{\overset{\_}{\kappa }^{2}+\overset{\_}{%
\tau }^{2}}}{\overset{\_}{\tau }},~\tau _{\overset{\_}{N}}=\frac{\overset{\_}%
{\kappa }^{2}}{\overset{\_}{\tau }(\overset{\_}{\kappa }^{2}+\overset{\_}{%
\tau }^{2})}(\frac{\overset{\_}{\tau }}{\overset{\_}{\kappa }})^{\prime },
\label{eqn26}
\end{equation}

\begin{equation}
\kappa _{\overset{\_}{B}}=\frac{\sqrt{\overset{\_}{\kappa }^{2}+\overset{\_}{%
\tau }^{2}}}{\overset{\_}{\tau }},~\tau _{\overset{\_}{B}}=\frac{\overset{\_}%
{\kappa }^{2}}{\overset{\_}{\tau }(\overset{\_}{\kappa }^{2}+\overset{\_}{%
\tau }^{2})}(\frac{\overset{\_}{\tau }}{\overset{\_}{\kappa }})^{\prime },
\label{eqn27}
\end{equation}%
respectively \cite{Uzunoglu}. Then, we can give the followings.

\begin{theorem}
The statements given below are equivalent.

i) The tangent inticatrix $(\overset{\_}{T})$ of $\beta $ is general helix.

ii) Osculating mate $\beta $ is slant helix.

iii) $\alpha $ is general helix.
\end{theorem}

\begin{proof}
Writing (\ref{eqn4}) into (\ref{eqn25}) and considering (\ref{eqn1}), it
follows $\frac{\tau _{\overset{\_}{T}}}{\kappa _{\overset{\_}{T}}}=\overset{%
\_}{\sigma }=\varepsilon _{1}\frac{\kappa }{\tau },$ which finishes the
proof.
\end{proof}

\begin{theorem}
The statements given below are equivalent.

i) The principal normal inticatrix $(\overset{\_}{N})$ of $\beta $ is
general helix.

ii) Osculating mate $\beta $ is $C-$slant helix.

iii) $\alpha $ is slant helix.
\end{theorem}

\begin{proof}
Writing (\ref{eqn4}) into (\ref{eqn26}) and considering (\ref{eqn1}) and (%
\ref{eqn1-1}), we obtain $\frac{\tau _{\overset{\_}{N}}}{\kappa _{\overset{\_%
}{N}}}=\frac{1}{\overset{\_}{\mu }}=-\sigma ,$ which gives the desired
results.
\end{proof}

\begin{theorem}
The statements given below are equivalent.

i) The binormal inticatrix $(\overset{\_}{B})$ of $\beta $ is general helix.

ii) Osculating mate $\beta $ is slant helix.

iii) $\alpha $ is general helix.
\end{theorem}

\begin{proof}
Writing (\ref{eqn4}) into (\ref{eqn27}) and considering (\ref{eqn1}), we
have $\frac{\tau _{\overset{\_}{B}}}{\kappa _{\overset{\_}{B}}}=-\overset{\_}%
{\sigma }=-\varepsilon _{1}\frac{\kappa }{\tau },$ which gives the desired
statements.
\end{proof}

\subsection{Osculating type(OT) Osculating Mates}

In this subsection we define osculating-type osculating mate (or $OT$%
-osculating mate) in $E^{3}$ and give the relationships between osculating
mates and $OT$-osculating mates. This section also gives a method to obtain
a rectifying curve.

Given a space curve $\alpha :I\rightarrow E^{3}$ with Frenet triangle $%
\left\{ T,N,B\right\} $ and curvatures $\kappa ,~\tau $. The vector $\overset%
{\sim }{D}=\frac{\tau }{\kappa }(s)T(s)+B(s)$ is named as modified Darboux
vector of $\alpha $ \cite{Izumiya1}. Let now the curve $\alpha $ be a Frenet
curve and the curve $\beta $ be an osculating mate of $\alpha $. The curve $%
\beta $ is called osculating-type osculating mate (or $OT$-osculating mate)
of $\alpha $, if the position vector of $\beta $ is always contained in the
osculating plane of $\alpha $.

Considering the definition of $(OT)$-osculating mate, we can write 
\begin{equation}
\beta (s)=m(s)T(s)+n(s)N(s),  \label{eqn28}
\end{equation}%
where $m(s),~n(s)$ are non-zero smooth functions of $s$. From (\ref{eqn3}), 
\begin{equation}
\left\{ 
\begin{array}{c}
T=\sin \left( \int \kappa ds\right) \overset{\_}{T}+\cos \left( \int \kappa
ds\right) \overset{\_}{B}, \\ 
N=\cos \left( \int \kappa ds\right) \overset{\_}{T}-\sin \left( \int \kappa
ds\right) \overset{\_}{B.}%
\end{array}%
\right.  \label{eqn29}
\end{equation}%
Writing (\ref{eqn29}) in (\ref{eqn28}) gives 
\begin{equation}
\left\{ 
\begin{array}{c}
\beta (s)=\left[ m\sin \left( \int \kappa ds\right) +n\cos \left( \int
\kappa ds\right) \right] \overset{\_}{T} \\ 
+\left[ m\cos \left( \int \kappa ds\right) -n\sin \left( \int \kappa
ds\right) \right] \overset{\_}{B.}%
\end{array}%
\right.  \label{eqn30}
\end{equation}%
Defining 
\begin{equation}
\left\{ 
\begin{array}{c}
\zeta (s)=m\sin \left( \int \kappa ds\right) +n\cos \left( \int \kappa
ds\right) , \\ 
\eta (s)=m\cos \left( \int \kappa ds\right) -n\sin \left( \int \kappa
ds\right) ,%
\end{array}%
\right.  \label{eqn31}
\end{equation}%
in (\ref{eqn30}) and differentiating the obtained equality gives%
\begin{equation}
\overset{\_}{T}=\zeta ^{\prime }\overset{\_}{T}+(\zeta \overset{\_}{\kappa }%
-\eta \overset{\_}{\tau })\overset{\_}{N}+\eta ^{\prime }\overset{\_}{B}.
\label{eqn32}
\end{equation}%
Hence, we get 
\begin{equation}
\eta =a=const,~\zeta =s+b=\frac{\overset{\_}{\tau }}{\overset{\_}{\kappa }}a,
\label{eqn33}
\end{equation}%
where $a,~b$ are non-zero constants. Considering (\ref{eqn33}), we obtain%
\begin{equation}
\beta (s)=a\left( \frac{\overset{\_}{\tau }}{\overset{\_}{\kappa }}\overset{%
\_}{T}+\overset{\_}{B}\right) (s)=a\overset{\sim }{\overset{\_}{D}}(s),
\label{eqn34}
\end{equation}%
where $\overset{\sim }{\overset{\_}{D}}$ is the modified Darboux vector of $%
\beta .$Then, the following theorem is obtained.

\begin{theorem}
Let $\beta $ be $OT$-osculating mate of $\alpha .$ Then,

i) $\beta $ is rectifiyng curve.

ii) The position vector $\beta $ and modified Darboux vector $\overset{\sim }%
{\overset{\_}{D}}$ of osculating mate $\beta $ are linearly dependent.
\end{theorem}

\bigskip

Considering (\ref{eqn31}), (\ref{eqn33}) and (\ref{eqn29}), the last theorem
gives a method to construct a rectifying curve by using osculating mate as
follows:

\begin{corollary}
The curve $\beta $ given by the parametrization 
\begin{equation}
\left\{ 
\begin{array}{c}
\beta (s)=\left[ (s+b)\sin \left( \int \kappa ds\right) +a\cos \left( \int
\kappa ds\right) \right] T(s) \\ 
+\left[ (s+b)\cos \left( \int \kappa ds\right) -a\sin \left( \int \kappa
ds\right) \right] N(s)%
\end{array}%
\right.  \label{eqn35}
\end{equation}%
is a rectifying curve and also osculating mate of $\alpha ,~$where $a,~b$
are non-zero constants.
\end{corollary}

\begin{example}
Let consider the spherical helix $\alpha $ in $E^{3}$ defined by 
\begin{equation*}
\alpha \left( t\right) =\left( \frac{1}{\sqrt{2}}\sin t,\cos t\cos (\sqrt{2}%
t)+\frac{1}{\sqrt{2}}\sin t\sin (\sqrt{2}t),-\cos t\sin (\sqrt{2}t)+\frac{1}{%
\sqrt{2}}\sin t\cos (\sqrt{2}t)\right) .
\end{equation*}%
(Fig. 1 (a)). The arc parameter of $\alpha $ is $s=\sin t$. The Frenet
apparatus of $\alpha $ are computed as follows,%
\begin{eqnarray*}
T\left( s\right) &=&\frac{\sqrt{2}}{2}\left( 1,-\sin (\sqrt{2}\arcsin
s),~-\cos (\sqrt{2}\arcsin s)\right) , \\
N\left( s\right) &=&\left( 0,-\cos (\sqrt{2}\arcsin s),~\sin (\sqrt{2}%
\arcsin s)\right) , \\
B\left( s\right) &=&-\frac{\sqrt{2}}{2}\left( 1,\sin (\sqrt{2}\arcsin
s),~\cos (\sqrt{2}\arcsin s)\right) , \\
\kappa &=&\frac{1}{\sqrt{1-s^{2}}},~~\ ~\tau =-\frac{1}{\sqrt{1-s^{2}}}.
\end{eqnarray*}%
From (\ref{eqn3})\ and (\ref{eqn8}), the osculating mate $\beta $ of $\alpha 
$ is obtained as 
\begin{equation*}
\beta (s)=\dint \left( sT(s)+\cos (\arcsin s)N(s)\right) ds=\left( \beta
_{1}\left( s\right) ,\beta _{2}\left( s\right) ,\beta _{3}\left( s\right)
\right) ,
\end{equation*}%
where%
\begin{eqnarray*}
\beta _{1}\left( s\right) &=&\frac{\sqrt{2}}{4}s^{2}+c_{1}, \\
\beta _{2}\left( s\right) &=&\dint \left( -\frac{\sqrt{2}}{2}s\sin (\sqrt{2}%
\arcsin s)-\cos (\sqrt{2}\arcsin s)\cos (\arcsin s)\right) ds, \\
\beta _{3}\left( s\right) &=&\dint \left( -\frac{\sqrt{2}}{2}s\cos (\sqrt{2}%
\arcsin s)+\cos \left( \arcsin s\right) \sin (\sqrt{2}\arcsin s)\right) ds,
\end{eqnarray*}%
where $c_{1}$ is integration constant. (Fig. 1(b)). From Theorem 14, the
osculating mate $\beta $ is a slant helix and its tangent indicatrix $%
\overset{\_}{T}$ is a general helix which is ploted in Fig. 2(a).
Furthermore, by choosing $a=b=\sqrt{2},$ from (\ref{eqn35}) an OT-osculating
mate of \ $\alpha $ is obtained easily which is also a rectifying curve and
ploted in Figure 2 (b).

\FRAME{itbpFU}{1.8343in}{1.8343in}{-0.0104in}{\Qcb{Fig. 1(a) Spherical helix 
$\protect\alpha $}}{\Qlb{Figure 1}}{fig1(a).jpg}{\special{language
"Scientific Word";type "GRAPHIC";maintain-aspect-ratio TRUE;display
"USEDEF";valid_file "F";width 1.8343in;height 1.8343in;depth
-0.0104in;original-width 4.1667in;original-height 4.1667in;cropleft
"0";croptop "1";cropright "1";cropbottom "0";filename
'Fig1(a).jpg';file-properties "XNPEU";}} \ \ \ \ \ \ \ \ \FRAME{itbpFU}{%
1.8343in}{1.8343in}{0in}{\Qcb{Fig. 1(b) Osculating mate $\protect\beta $}}{%
\Qlb{Figure 2}}{fig1(b).jpg}{\special{language "Scientific Word";type
"GRAPHIC";maintain-aspect-ratio TRUE;display "USEDEF";valid_file "F";width
1.8343in;height 1.8343in;depth 0in;original-width 4.1667in;original-height
4.1667in;cropleft "0";croptop "1";cropright "1";cropbottom "0";filename
'Fig1(b).jpg';file-properties "XNPEU";}}

\FRAME{itbpFU}{1.8343in}{1.8343in}{0in}{\Qcb{Fig. 2(a) Tangent indicatix $%
\protect\overset{\_}{T}$}}{}{fig1(c).jpg}{\special{language "Scientific
Word";type "GRAPHIC";maintain-aspect-ratio TRUE;display "USEDEF";valid_file
"F";width 1.8343in;height 1.8343in;depth 0in;original-width
4.1667in;original-height 4.1667in;cropleft "0";croptop "1";cropright
"1";cropbottom "0";filename 'Fig1(c).jpg';file-properties "XNPEU";}} \ \ \ \
\ \ \ \ \ \FRAME{itbpFU}{1.8343in}{1.8343in}{0in}{\Qcb{Fig. 2(b)
OT-osculating mate of $\protect\alpha $}}{}{fig2(b).jpg}{\special{language
"Scientific Word";type "GRAPHIC";maintain-aspect-ratio TRUE;display
"USEDEF";valid_file "F";width 1.8343in;height 1.8343in;depth
0in;original-width 4.1667in;original-height 4.1667in;cropleft "0";croptop
"1";cropright "1";cropbottom "0";filename 'Fig2(b).jpg';file-properties
"XNPEU";}} \ 
\end{example}

\section{Conclusions}

A new type of associated curves is introduced and named as osculating mate.
The relations between a Frenet curve and its osculating mate are obtained.
The obtained results allow to construct a slant helix, a $C$-slant helix, a
spherical helix and a rectifying curve by considering osculating mate of a
Frenet curve.

\bigskip

\section{\protect\bigskip Compliance with Ethical Standards}

Funding: Not aplicable. (There is no funding)

Conflict of Interest: The authors wish to confirm that there are no known
conflicts of interest associated with this publication and there has been no
significant financial support for this work that could have influenced its
outcome.

Data Availability Statement: The authors wish to confirm that this
manuscript has no associated data.

Ethical approval: This article does not contain any studies with animals
performed by any of the authors.

Ethical approval: This article does not contain any studies with human
participants or animals performed by any of the authors.

\end{document}